\magnification=1200
\input mssymb
\font\kapit=cmcsc10
\def\n{n<\omega}

\font\eightslss = cmssi8
\def \qed{\hskip10dd {\eightslss QED}}
\def\GG{\hbox{\bf G}}

\def\dm{\hbox{dm}}
\def\PP{\hbox{\bf P}}
\def\RR{\hbox{\bf R}}
\def\QQ{\hbox{\bf Q}}

\def\T{T^{(n)}_{ni}}
\def\t{T^{(n)}}
\def\q{q^{(n)}}
\def\o{\omega}

\centerline{On closed $P$-sets with $ccc$ in the space $\omega^*$}

\bigskip
\centerline{R.Frankiewicz(Warsaw) S.Shelah(Jerusalem) P.Zbierski(Warsaw)}
\bigskip
Abstract. It is proved that -- consistently -- there can be no $ccc$
closed $P$-sets in the remainder space $\omega^*$.

\bigskip

In this paper we show how to construct a model of set theory in which
there are no $P$-sets satisfying $ccc$ (countable antichain condition)
in the ultrafilter space $\omega^*=\beta\omega|\setminus\omega$.
The problem of the existence of such sets (which are generalizations
of $P$-points) was known since some time and occurred explicitly in
vM-R|. In the proof we follow the construction from S| of a model
in which there are no $P$-points. A particular case of $P$-sets,
which are supports of approximative measures has been settled in M|,
where  the author shows that there can be no such measures on
$P(\omega)/fin$. (Under $CH$, e.g. the Gleason space $\GG(2^\omega)$ of
the Cantor space is a ccc $P$-set in $\o^*$ which carries no
approximative measure).

\bigskip
Sec.1. Closed $P$-sets in the space $\o^*$ can be identified with
$P$-filters $F$ on $\o$. Thus, the dual ideal   $I=\{\omega\setminus
A:\ A\in F\}$ has the property:

$$\eqalign{&\hbox{If $A_n\in I$,  for $n\in\omega$, then
there is an $A\in I$}\cr
&\hbox{such that
$A_n\subseteq_* A$, for each $n\in\omega$.}\cr}\leqno(1.1)$$

Further, the countable chain condition imposed upon $F$ implies that
$I$ is fat in the following sense (see F-Z|):

$$\eqalign{&\hbox{if $A_n\in I$, for $n\in\omega$ and
$\lim_n\min A_n=\infty$,}\cr
&\hbox{then
there is an infinite $Z\subseteq\omega$ such that $\bigcup_{n\in
Z}A_n\in I$.}\cr}\leqno(1.2)$$

Indeed, let $e_n=A_n\setminus A$, where $A\in I$ is as in (1.1).
Since $\min A_n$ are arbitrarily large, we can find an infinite set
$Y\subseteq\omega$ such that the family $\{e_n: n\in Y\}$ is
disjoint. If $\{Y_\alpha:\ \alpha<c\}$ is an almost disjoint family of
subsets of $Y$, then the unions

$$S_\alpha=\bigcup\{e_n:\ n\in Y_\alpha\},\ \ \ \alpha<c$$

\noindent are almost disjoint and hence the closures $S^*_\alpha$ in
the space $\o^*$ are disjoint. By $ccc$ we have

$$S_\alpha^*\cap\bigcap\{B^*:\ B\in F\}=\emptyset,$$

\noindent for some $\alpha$ and consequently $S_\alpha\in I$. It follows that th
e
union

$$\bigcup_{n\in Y_\alpha}A_n=\bigcup_{n\in Y_\alpha}(A_n\cap
A)\cup\bigcup_{n\in Y_\alpha}(A_n\setminus A)$$

\noindent is in $I$ as a subset of $S_\alpha\cup A$.

Let us fix a given ccc $P$-filter $F$ and its dual $I$. We shall
define a forcing $\PP=\PP(F)$.

A partial ordering $(T,\le_T)$, where $T\subseteq\omega$, will be
called a tree, if for each $i\in T$ the set of predecessors $\{j\in
T:\ j\le_T i\}$ is linearly ordered and

$$\hbox{$i\le_T j$ implies $i\le j$, for all $i,j\in T$.}$$

We define a partial ordering for trees

$$\eqalign{&\hbox{$T\le_t S$ iff $(S,\le_S)$ is a subordering of
$(T,\le_T)$}\cr
&\hbox{and each
branch of $T$ contains cofinally a (unique) branch of $S$.}\cr}$$

There is a tree $T_0$ such that $T_0\in I$ and $T_0$ is order
isomorphic to the full binary tree of height $\omega$.

Deleting the numbers $\le n$ from $T_0$ we obtain a subtree denoted by
$T_0^{(n)}$ (we have $T_0^{(n)}\le_t T_0^{(m)}$, for $n\le m$). Let
${\cal T}$ consist of all the trees $T\in I$ such that

$$\hbox{$T\le_t T_0^{(n)}$, for some $n\in \omega$.}$$

Note that each tree $T\in {\cal T}$ has finitely many roots.

\medskip
{\kapit Definition.} {\it Elements of the forcing $\PP$ are of the form
$p=<T_p,f_p>$, where $T_p\in\cal T$ and $f_p:\
T_p\longrightarrow\{0,1\}$. The ordering of $\PP$ is defined thus

$$\hbox{$p\le q$ iff $T_p\le_t T_q$ and $f_p\supseteq f_q$.}$$}
\medskip
Let $\{b_\alpha:\ \alpha<c\}$ be a fixed enumeration of all the
branches of $T_0$ in $V$. For a generic $G\subseteq\PP$ let
$T_G=\bigcup_{p\in G} T_p$ and $f_G=\bigcup_{p\in G}f_p$.

Each branch $B$ of $T_G$ contains cofinally a unique $b_\alpha$. Let us
write $B=B_\alpha$ and define

$$X_\alpha=\{i\in\omega:\ i\in B_\alpha\ \hbox{and}\ f_G(i)=1\}$$

Since $T_p\in I$, for any $p\in \PP$, hence $\o\setminus T_p\cap A$ is infinite,
for each $A\in F$. It follows that the sets

$$D_{n\varepsilon}^{A\alpha}=\{p\in\PP:\ \exists i>ni\in b^p_\alpha\
\hbox{and}\ f_p(i)=\varepsilon|\}$$

\noindent are dense, for each $A\in F$, $n\in\omega$, $\alpha<c$ and
$\varepsilon=0,1$ (here $b_\alpha^p$ denotes the branch of $T_p$
extending $b_\alpha$).

Thus, $\PP$ adds uncountably many almost disjoint Gregorieff-like sets.

\bigskip
Sec.2. Let $\QQ=\QQ(F)$ be a countable product of $\PP=\PP(F)$. Thus
the elements $q\in\QQ$ can be written in the form

$$\hbox{$q=<f_0^q, f_1^q,\ldots>$, where
$<\hbox{dm}(f_i^q),f_i^q>\in\PP$, for each $i<\o$}.$$

By $q^{(n)}$ we denote the condition $<g_i:\ i<\omega>$ where

$$g_i=\cases{f_i^q\mid\dm(f_i^q)^{(n)},\ &for$\ i<n$\cr
f_i^q&for$\ i\ge n$\cr}$$

Here $T^{(n)}$ is a tree obtained from $T$ by deleting the numbers $\le n$.

\medskip
{\kapit Lemma 2.1.}{\it For each decreasing sequence $p_0\ge
p_1\ge\ldots$ there
is a $q$ and an infinite $Z\subseteq\omega$ such that

$$\hbox{$q\le p_n^{(n)}$, for each $n\in Z$}.$$}

\smallskip
{\kapit Proof.} Let $T_{ni}=\dm(f_i^{p_n})$, where $p_n=<f_i^{p_n}:\
i<\omega>$. Since $\min\T\ge n$, we may use (1.2) to define by
induction a descending sequence $Z_0\supseteq Z_1\supseteq\ldots$ of
infinite subsets of $\omega$ such that

$$\hbox{$\bigcup_{n\in Z_i}\T$ is in $I$, for each $i<\omega$}.$$

There is an infinite $Z\subseteq\omega$, such that $Z\subseteq_*Z_i$,
for each $i<\omega$. Define

$$T_i=T_{ii}\cup\bigcup_{n\in Z}\T$$

\noindent and

$$f_i^q=f_i^{p_i}\cup\bigcup_{n\in Z}f_i^{p_n}\mid\T.$$

Then, $\dm(f^q_i)=T_i$ and $q=<f_i^q:\ i<\omega>$ is as required.\qed

\medskip
For $q\in\QQ$ and $n\in\omega$ let $S(q,n)$ be the set of all
sequences $s=<s_0,\ldots,s_{n-1}>$ satisfying the following properties

\smallskip
\item{1.} $s_0,\ldots,s_{n-1}$ are finite zero-one functions.

\item{2.} The domains
$t_0=\dm(s_0),\ldots,t_{n-1}=\dm(s_{n-1} )$ are finite
trees such that

$$t_0\cap T^{(n)}_0=\ldots=t_{n-1}\cap T^{(n)}_{n-1}=\emptyset,$$
where $T_0=\dm(f_0^q),\ldots,T_{n-1}=\dm(f_{n-1}^q)$

\item{3.} Ordered sums
$t_0\oplus T^{(n)}_0,\ldots,t_{n-1}\oplus \t_{n-1}$ are
trees belonging to $\cal T$.

\smallskip
Note that from the definition of $\cal T$ it follows that $S(q,n)$ is
always finite. Let us denote

$$s*q^{(n)}=<s_0\cup f_0^q,\ldots,s_{n-1}\cup f_{n-1}^q,f_n^q,\ldots>$$

\noindent for $q,n,s$ as above. Obviously, we have

 $$\hbox{the set $\{s*q^{(n)}:\ s\in S(q,n)\}$ is
predense below $\q$}\leqno(2.2)$$

\noindent(i.e. the boolean sum $\sum_{s\in S(q,n)}s*\q=\q$).

Now, we obtain easily an analogue of VI, 4.5 in S|.

$$\eqalign{&\hbox{For arbitrary $p\in\QQ$, $n<\omega$ and
$\tau\in V^{(Q)}$}\cr
&\hbox{such
that $\QQ\Vdash``\tau\ \hbox{is an ordinal"}$ there is}\cr
&\hbox{a $q\le p$ and
ordinals $\{\alpha(s):\ s\in S(p,n)\}$}\cr
&\hbox{so that}\cr
&\q\Vdash``\bigvee_s\tau=\alpha(s)"\cr}\leqno(2.3)$$

Indeed, if $S(q,n)=\{s^0,\ldots,s^{m-1}\}$, then we define
inductively conditions $p_0,\ldots,p_{m}$ so that $p_0=p$ and
$p_{k+1}\le s^k*p_k^{(n)}$ is such that

$$\hbox{$p_{k+1}\Vdash``\tau=\alpha"$, for some ordinal
$\alpha=\alpha(s^k)$}$$

Now, $q=s*p^{(n)}_{m}$, where $s$ is such that $p=s*p^{(k)}$ (we
may assume $s\in S(p,n))$, satisfies (2.3).

\medskip
\noindent 2.4. {\kapit Theorem.}{\it $\QQ$ is $\alpha$-proper, for
every $\alpha<\omega_1$,
and has the strong $PP$-property.}

\smallskip
{\kapit Proof.} Let countable ${N}\prec H(\kappa)$ for sufficiently large
$\kappa$, be such that $\QQ\in{N}$ and suppose that
$p\in\QQ\cap{N}$. To prove that $\QQ$ is proper we have to find
a $q\le p$, which is $N$-generic. Let $\{\tau_n:\ n<\omega\}$ be
an enumeration of all the $\QQ$-names for ordinals, such that
$\tau_n\in{N}$, for $n<\omega$. Using (2.3) we define
inductively a sequence $p_0=p\ge p_1\ge\ldots$ and ordinals
$\alpha(n,s)$ so that

$$p_n^{(n)}\Vdash``\bigwedge_{i\le n}\bigvee_s\tau_i=\alpha(n,s)"\
\hbox{for each $n<\omega$}$$

\noindent (i.e. in the n-th step we apply (2.3) for all names
$\tau_0,\ldots,\tau_n$). Note that the $p's$ and $\alpha's$ can be
found in $N$, since ${N}\prec H(\kappa)$. By Lemma 2.1
there is a $q$ and an infinite $Z\subseteq\omega$ such that

$$\hbox{$q\le p_n^{(n)}$, for each $n\in Z$}.$$

Hence also $q^{(m)}\le p_n^{(n)}$ holds for arbitrarily large $n$ and
all $m<\omega$ and thus $$q^{(m)}\Vdash``\tau_n\in{N}",$$ for all
$n,m<\omega$.

By III, 2.6 of S|, each $q^{(m)}$ is ${N}$-generic.

To see that $\QQ$ is $\alpha$-proper let $<N_\xi:\ \xi\le\alpha>$
be a continuous sequence of elementary countable submodels of
$H(\kappa)$ such that $\QQ\in N_0$ and

$$<N_\xi:\ \xi\le\eta>\in N_{\eta+1},\ \hbox{for each}\ \eta<\alpha.$$

Assume that $\QQ$ is $\beta$-proper, for each $\beta<\alpha$ and let
$q_0\in\QQ\cap{N}_0$. If $\alpha=\beta+1$, we have a $q\le q_0$
which is $N_\xi$-generic, for each $\xi\le\beta$ and we may assume
that the $\q$ have the same property, for all $n<\omega$. since
$N_\alpha\prec H(\kappa)$ and all the parameters are in $N_\alpha$,
such a $q$ can be found in $N_\alpha$ and as above we construct a
$q_\alpha\le q$ which is $N_\alpha$-generic and so are the
$q_\alpha^{(n)}$, for $n<\omega$. Thus, $q_\alpha$ and all the
$q_\alpha^{(n)}$ are $N_\xi$-generic for all $\xi\le\alpha$. If
$\alpha$ is a limit ordinal, we fix an increasing sequence $<\xi_n:\
n<\omega>$ such that $\alpha=\sup_{\n}\xi_n$ and by the inductive
hypothesis there is a sequence $q_0\ge q_{\xi_0}\ge
q_{\xi_1}\ge\ldots$ such that, for each $\n$, $q_{\xi_n}$ is
$N_\xi$-generic, for each $\xi\le\xi_n$ and $q_{\xi_n}\in
N_{\xi_n+1}$ and that $q_{\xi_n}^{(m)}$ have the same property for
each $m<\omega$. By Lemma 2.1 there is a $q\in\QQ$ such that $q\le
q_{\xi_n}^{(n)}$, for infinitely many $\n$. Thus, $q\le q_0$ and $q$
is $N_\xi$-generic for each $\xi<\alpha$ and hence also for each
$\xi\le\alpha$.

Finally, to prove the $PP$-property let
$h:\omega\longrightarrow\omega$ diverge to infinity and suppose that
$p\Vdash ``f:\ \omega\longrightarrow\omega"$. Define

$$k_n=\min\{i:\ h(i)>2^n\cdot\hbox{card}S(p,n)\},\ \hbox{for}\ \n$$

\noindent and, using (2.3), define inductively the sequence $p=p_0\ge
p_1\ge\ldots$ such that

$$p_n^{(n)}\Vdash``\bigwedge_{i<k_n}\bigvee_{s\in
S(p_i,i)}f(i)=\alpha(s,i)"\ \hbox{for each $\n$ and some integers
$\alpha(s,i)<\omega$}.$$

Let $T$ be the tree built up of integers

$$\{\alpha(s,i):\ i<\omega\ \hbox{and}\ s\in S(p_i,i)\}$$.

\noindent If $q\le p_n^{(n)}$, for infinitely many $n$, then we have
$q\Vdash``f\in \hbox{Lim}\ T"$ and $T\cap\omega^{k_n}$ has less
elements than $h(k_n)$, for all $\n$, which finishes the proof.\qed

\medskip
The last point to be discussed is how does $\QQ=\QQ(F)$ act in the
course of iteration.

\medskip
\noindent 2.5.{\kapit Lemma.}{\it If $\RR$ is $\o^\o$-bounding (i.e.
the set of old
functions $:\o\longrightarrow\o$ dominates) and $\QQ(F)$ is a
complete subforcing of $\RR$, then in $V^{(\RR)}$ the filter $F$ cannot be
extended to a $ccc$ $P$-filter.}

\smallskip
{\kapit Proof.} Let $X^n_\alpha$ be the $\alpha$-th set added by $n$-th factor
of the product $\QQ=\PP^\omega$. Suppose that for some $r\in\RR$ and
a $ccc$ $P$-filter $E\in V^{(\RR)}$ we have

$$r\Vdash``F\subseteq E"$$

Note that for each $n<\omega$, the relation $X^n_\alpha\in E$ hold
for at most countably many $\alpha$'s, since $E$ is $ccc$. Hence,
there is an $\alpha$ such that for all $\n$ we have $\o\setminus
X_\alpha^n\in E$ and, since $E$ is a $P$-filter, there is an $A\in E$
and a function $g$, so that $A\subseteq\bigcap_{\n}(\o\setminus
X_\alpha^n)\cup0,g(n))$ i.e. for some $r_1\le r$ we have

$$r_1\Vdash``\bigcap_{\n}(\o\setminus X_\alpha^n)\cup0,g(n))\in E"\leqno(2.6)$$

Since $\RR$ is $\o^\o$-bounding we may assume that $g\in V$. By the
assumption, $\QQ$ is a complete subforcing of $\RR$ and hence there
is a $q\in\QQ$ such that $r$ is compatible with each $q'\le q$.

On the other hand, since $T_n=\dm(f_n^q)\in{\cal T}$, there is a set
$B\in I$ and an increasing sequence $a_0<a_1<\ldots$ such that
$T_n\setminus0,a_n)\subseteq B$, $g(n)<a_n$ and
$a_n,a_{n+1})\setminus B\not=\emptyset$, for each $\n$. Define
$q'\le q$ as follows. For a given $n$ extend $T_n$ by adding elements
of $a_n,a_{n+1})\setminus B$ on the $\alpha$-th branch $b_\alpha^q$
and put $f_n^{q'}(i)=1$, for each $i\ina_n,a_{n+1})\setminus B$.
Obviously, we have

$$q'\Vdash``(\o\setminus
X_\alpha^n)\cup0,g(n))\capa_n,a_{n+1})\setminus B=\emptyset",\
\hbox{for each $n$}$$

\noindent and hence

$$q'\Vdash``\bigcap_{\n}(\o\setminus X_\alpha^n)\cup0,g(n))\setminus
B\cap\bigcup_{\n}a_n,a_{n+1})=\emptyset"$$

Consequently $q'\Vdash``\bigcap_{\n}(\o\setminus
X_\alpha^n)\cup0,g(n))\subseteq_* B"$, which contradicts (2.6).\qed

\medskip
The rest of the proof is routine. Beginning with a model $V$ of
$2^\o=\o_1$ and $2^{\o_1}=\o_2$ we iterate with countable supports
the forcings $\QQ(F)$, for all $ccc$ $P$-filters $F$ booked at each
stage $\alpha<\o_2$ of the iteration. From S|, V.4 we know that the
resulting forcing $\RR$ (obtained after $\o_2$ stages) is proper and
$\o^\o$-bounding. Hence, in $VG|$ there are no $ccc$ $P$-sets.

\bigskip
\centerline{References}

F-Z| R.Frankiewicz, P.Zbierski ``Strongly discrete subsets in $\o^*$"
Fund. Math. 129 (1988) pp 173-180

M|A.H.Mekler ``Finitely additive measures on N and the additive
property" \break Proc.AMS Vol.92 No3 Nov.1984 pp 439-444

vM-R| J.vanMill, G.M.Reed editors ``Open problems in topology",
North Holland Elsevier Science Publishers B.V 1990

S| S.Shelah ``Proper forcing" Lecture Notes in Mathematics 940
Springer Verlag

\bigskip

{\kapit Ryszard Frankiewicz} Institute of Mathematics, Polish Academy
of Sciences, Warsaw, Poland

{\kapit Saharon Shelah} Institute of Mathematics, The Hebrew
University, Jerusalem, Israel

{\kapit Pawe\l Zbierski} Institute of Mathematics,  Warsaw
University, Warsaw, Poland\bye